\documentclass[twoside,12pt]{article}
\usepackage[]{amsmath,amsthm,amssymb}
\usepackage[latin1]{inputenc}

\begin{document}
\title{{\Large\bf  On the Yor integral and a system of polynomials related to the Kontorovich-Lebedev transform}}

\author{Semyon   YAKUBOVICH}
\maketitle

\markboth{\rm \centerline{ Semyon YAKUBOVICH}}{}
\markright{\rm \centerline{THE YOR INTEGRAL AND A SYSTEM OF
POLYNOMIALS}}

\begin{abstract} {In this paper we establish different representations of the so-called Yor integral, which is one of the
key ingredient in mathematical finance, in particular, to compute
normalized prices of Asian options. We show, that the Yor integral
is related with the Kontorovich-Lebedev transform. Also we discuss
its relationship with a system of polynomials recently introduced by
the author. We  derive  new important properties of these polynomials, including upper bounds,  an exact asymptotic behavior for large values of their degree and explicit formula of coefficients.}

\end{abstract}
\vspace{4mm}

{\bf Keywords}: {\it Yor integral,  Kontorovich-Lebedev transform,
Parseval equality, modified Bessel functions, Asian options,
Hartman-Watson distribution,  Brownian motion}

{\bf AMS subject classification}:  44A15,   60J65, 60E99

\vspace{4mm}

\section {Introduction and preliminary results}

In 1980 Yor [14] (see also in [3], [4], [15]) expressed the density
of the Hartman -Watson distribution [2], which is related to the
pricing of Asian options in mathematical finance, in the form of
elementary integral, involving exponential functions
$$F_t(r)= \frac{ e^{\pi^2/2t}}{ \sqrt {2\pi^3t}}\int_0^\infty
\exp\left(- {y^2\over 2t}\right)\exp\left(-r\cosh y\right)\sinh y
\sin\left({\pi y\over t}\right) dy,\  r,t > 0.\eqno(1.1)$$
However, despite of the importance of integral (1.1) in applications
and a necessity to calculate it in the closed form, this task  is
quite difficult. Nevertheless, we can perhaps stimulate the
corresponding numerical calculations of the integral and finding its
asymptotic behavior for small values of $t$, representing (1.1) in a
different form. Precisely, we will show the relation between
$F_t(r)$ and the Kontorovich-Lebedev transformation [9], [10], [11]
$$(Gf)(\tau)=\int_0^\infty K_{i\tau}(r)f(r) dr, \ \tau \in \mathbb{R}_+,\eqno(1.2)$$
where $K_{i\tau}(r)$ is the modified Bessel function of the pure
imaginary index $i\tau$ [1]. As it is known,  operator (1.2) is
bounded
$$G: L_2\left({\mathbb {R}}_+;  rdr \right)
\leftrightarrow L_2({\mathbb {R}}_+;\tau\sinh\pi\tau d\tau)$$
and the integral in (1.2) converges with respect to the norm in
$L_2({\mathbb {R}}_+;\tau\sinh\pi\tau d\tau)$. Moreover, it forms an
isometric isomorphism between these Hilbert spaces and the Parseval
identity holds
$$\int_0^\infty \tau\sinh\pi\tau |(Gf)(\tau)|^2d\tau= {\pi^2\over 2}
\int_0^\infty |f(r)|^2 rdr .\eqno(1.3)$$
Reciprocally,  the inversion formula
$$f(r)={2\over r \pi^2}\int_{0}^\infty \tau\sinh\pi\tau K_{i\tau}(r)(Gf)(\tau) d\tau,\eqno(1.4)$$
takes place, where the convergence of the integral (1.4) is
understood with respect to the norm of the space $L_2\left({\mathbb
{R}}_+;  rdr \right)$.

The modified Bessel function $K_{i\tau}(x)$ is an eigenfunction of
the following second order differential operator
$${\cal A}_x \equiv x^2- x{d\over dx} x {d\over dx},\eqno(1.5)$$
i.e. we have
$${\cal A}_x \ K_{i\tau}(x)= \tau^2 K_{i\tau}(x).\eqno(1.6)$$
It has the asymptotic behavior (cf. [1] relations (9.6.8), (9.6.9),
(9.7.2))
$$ K_\nu(z) = \left( \frac{\pi}{2z} \right)^{1/2} e^{-z} [1+
O(1/z)], \qquad z \to \infty,\eqno(1.7)$$
and near the origin
$$ K_\nu(z) = O\left ( z^{-|{\rm Re}\nu|}\right), \ z \to 0,\eqno(1.8)$$
$$K_0(z) = -\log z + O(1), \ z \to 0. \eqno(1.9)$$
Moreover it can be defined by the  following integral
representations [9, (6-1-2)],  [8], Vol. I, relation (2.4.18.4)
$$K_{\nu}(x)=\int_0^\infty e^{-x\cosh u}\cosh \nu u du,
x>0, \eqno(1.10)$$
$$K_{\nu}(x)={1\over 2} \left(x\over 2\right)^\nu\int_0^\infty
e^{-t-{x^2\over 4t}}t^{-\nu-1}dt, \ x>0.\eqno(1.11)$$
The convolution operator for the Kontorovich-Lebedev transform is
defined as follows [10, 11]
$$(f*h)(x)\equiv (f(x)* g(x)) ={1\over 2x}\int_0^\infty\int_0^\infty e^{-{1\over
2}\left(x{u^2+y^2\over uy}+{yu\over x}\right)}f(u)h(y)dudy, \
x>0.\eqno(1.12)$$
It is well defined in the Banach ring $L^\alpha({\mathbb{R}}_+)
\equiv L_{1}({\mathbb{R}}_+;K_{\alpha}(x)dx), \alpha \in
{\mathbb{R}}$, i.e. the space of all summable functions $f:
{\mathbb{R}}_+ \to {\mathbb{C}}$ with respect to the measure
$K_{\alpha}(x)dx$ for which
 $$||f||_{L^\alpha({\mathbb{R}}_+)}=\int_0^\infty |f(x)|K_{\alpha} (x)dx\eqno(1.13)$$
is finite. The following embeddings take place
$$L^{\alpha}({\mathbb{R}}_+)\equiv L^{-\alpha}({\mathbb{R}}_+), \
L^{\alpha}({\mathbb{R}}_+) \subseteq L^{\beta}({\mathbb R}_+), \
|\alpha|\ge |\beta| \ge 0,  \alpha, \beta  \in {\mathbb{R}},$$
$$L^{\alpha}({\mathbb{R}}) \supset L_p({\mathbb{R}}_+; xdx),  \  2 < p \le \infty, \ |\alpha|< 1-
{2\over p},$$
where $L_p({\mathbb{R}}_+; xdx)$ is a weighted Banach space with the
norm
$$||f||_{L_p({\mathbb{R}}_+; xdx)}=\left(\int_0^\infty |f(x)|^p
xdx\right)^{1/p},  \  1 \le  p< \infty,$$
$$||f||_{L_\infty({\mathbb{R}}_+; xdx)}= \hbox{ess sup}_{x \in {\mathbb{R}}_+}|f(x)|.$$
The factorization  property is true for the convolution (1.12) in
terms of the Kontorovich-Lebedev transform (1.2)  in the space
$L^\alpha({\mathbb{R}}_+)$, namely
 $$(G[f*h])(\tau)=(Gf)(\tau)(Gh)(\tau),\  \tau \in {\mathbb R}_+.\eqno(1.14)$$
This property is based on the Macdonald formula [1]
$$K_\nu(x)K_\nu(y)= {1\over 2}\int_0^\infty e^{-{1\over
2}\left(t{x^2+y^2\over xy}+{xy\over t}\right)}K_\nu(t){dt\over
t}.\eqno(1.15)$$

\section{A system of polynomials}

In this section we will provide  a useful information about a system
of polynomials, which  is  related to the Kontorovich-Lebedev
transformation (1.2) and  was studied by the author in [12]. In
particular, we will derive new properties of these polynomials,
including an upper bound, a series representation by the index of
the modified Bessel functions of the third kind and an explicit
formula expressing  their coefficients. In fact, as it is proved in
[12], the following functions
$$p_n(x)= (-1)^n e^{x}{\mathcal A}_x^n\  e^{-x}, \ n \in
\mathbb{N}_0,\eqno(2.1)$$
where ${\mathcal A}_x^n$ is the $n$-th iteration of the differential
operator (1.5), are $n$-th degree polynomials, which have the
following integral representation
$$p_n(x)= {2(-1)^n\over \pi} e^{x}\int_{0}^\infty \tau^{2n}
K_{i\tau}(x) d\tau,\ x >0, \ n \in \mathbb{N}.\eqno(2.2) $$
Hence calling the reciprocal formula (1.2), we obtain the
Kontorovich-Lebedev transform of $p_n(x)$

$$\int_{0}^\infty  K_{i\tau}(x) e^{-x} p_n(x){dx\over x}=
(-1)^n \  {\pi \tau^{2n-1}\over \sinh(\pi\tau)}.\eqno(2.3)$$
The system $p_n$ satisfies the differential recurrence relation of
the form
$$p_{n+1}(x)= x^2p_n^{\prime\prime}(x) + x(1-2x)p_n^\prime(x)- xp_n(x), \ n= 0,1,2,\dots \ .$$
In particular, we derive
$$p_0(x)=1, \ p_1(x)= -x,\ p_2(x)= 3x^2-x, \quad p_3(x)= -15x^3+ 15x^2-x.$$
The generating function of these polynomials is given by the series
$$e^{- 2x\sinh^2( t/2)} = \sum_{n=0}^\infty {p_n(x)\over (2n)!}
\ t^{2n}.$$
Letting $x=0$ in the latter equation, we find
$$p_n(0)= 0, \quad n=1,2,\dots.$$
The leading coefficient $a_{n,n}$ of these polynomials can be calculated by the formula
$$a_{n,n}= (-1)^n (2n-1)!! =  (-1)^n 1\cdot 3\cdot 5\dots \  \cdot (2n-1), \ n \in \mathbb{N}.\eqno(2.4)$$

The following lemma proves an upper bound for the system $p_n(x)$.
Indeed, we have

{\bf Lemma 1}. {\it Let $ x >0,\ n \in \mathbb{N}, \ \varepsilon \in
\left(1- {\sqrt 3\over 2},1\right]$ and

$$\alpha \in \left(0, \quad 2\hbox{arccos}\ \frac{\sqrt{1+4(1-\varepsilon)^2}}{2}\ \right].\eqno(2.5)$$
Then }
$$|p_n(x)| \le \sqrt{\frac{(2^{4n}-1)(4n)!\sin(\alpha/2)}{\pi n
\ \alpha^{4n} (2^{4n-2}+ 6/\pi^2 -1)}}\ {e^{\varepsilon x}\over
2}.\eqno(2.6)$$

\begin{proof} In fact, taking representation (2.2), we apply the Schwarz inequality together with integral
formula for Bernoulli numbers $B_{4n}, n = 1,2,\dots$ \ [1]
$$\int_{0}^\infty {\tau^{4n-1}\over \sinh\tau }d\tau= - B_{4n}
\frac{(2^{4n}-1)\pi^{4n}}{4n},\eqno(2.7)$$
and relation (2.16.51.8) in [8], Vol. II. Then    we deduce
$$|p_n(x)| \le {2\over \pi} e^{x}\int_{0}^\infty |K_{i\tau}(x)|\tau^{2n}\ d\tau \le
{2\over \pi} e^{x}\left(\int_{0}^\infty \tau\sinh \left(\alpha
\tau\right)K_{i\tau}^2(x) \ d\tau\right)^{1/2}$$
$$\times \left(\int_{0}^\infty {\tau^{4n-1}\over
\sinh(\alpha\tau)}d\tau\right)^{1/2}= e^x \ \left({\pi\over
\alpha}\right)^{2n} \sqrt{\frac{(1-2^{4n})\sin(\alpha/2)B_{4n}}{2\pi n}}$$
$$\times \left(x K_1(2x\cos(\alpha/2))\right)^{1/2}, \ x > 0, \alpha \in (0, \pi).$$
But in the meantime via (1.10) it is not difficult to find that
$$\left(x K_1(2x\cos(\alpha/2))\right)^{1/2}= \left(x \int_0^\infty e^{-2x\cos(\alpha/2)\cosh u} \cosh u du
\right)^{1/2}$$
$$=\left(x \int_0^\infty e^{-2x\cos(\alpha/2)\sqrt {v^2+1}} dv
\right)^{1/2}\le e^{\varepsilon x- x} \left(x \int_0^\infty e^{-x v}
dv\right)^{1/2} = e^{\varepsilon x- x}$$
when $\varepsilon \in \left(1- {\sqrt 3\over 2},1\right]$ and
$\alpha$ satisfies condition (2.5). Indeed, in this case
$$2\cos(\alpha/2)\sqrt {v^2+1} \ge v + 2 (1-\varepsilon)$$
for any $v \ge 0$ since $2\cos(\alpha/2) \ge \sqrt{1+
4(1-\varepsilon)^2}.$ Moreover, employing  a sharp upper bound for
the Bernoulli numbers $-B_{4n}$ (see in [2])
$$-B_{4n}\le  {(4n)!\over 2\pi^{4n} (2^{4n-2}+ 6/\pi^2 - 1)},$$
we combine  with (2.7) to arrive at inequality (2.6) and complete
the proof of Lemma 1.
\end{proof}

{\bf Remark 1}. As we observe form the Stirling asymptotic formula
for factorials [1],  inequality (2.6) guarantees the absolute and
uniform convergence of  series (2.4) for the generating function on
any compact set of $x  \in [x_0, X_0] \subset \mathbb{R}_+$ and $t$
from the interval $|t|  \le t_0 < 2/ \alpha$.

{\bf Theorem 1.} {\it The system $p_n(x), x >0, \ n \in \mathbb{N}$
can be expressed in terms of the absolutely convergent series
$$p_n(x)= 2e^{x}\sum_{m=1}^\infty (-1)^{m}m^{2n}I_m(x),\eqno(2.8)$$
where $I_\nu(z)$ is the modified Bessel function of the third kind
[5].  Moreover, we have an explicit formula for these polynomials
$$p_n(x)= \sum_{k=1}^n a_{k,n} x^k,\eqno(2.9)$$
where the coefficients $a_{k,n}$ are given by
$$a_{k,n}= {1\over k!}\sum_{r=0}^{k} {(-1)^r\over 2^r} \binom{k}{r}\sum_{j=0}^{k-r} {(-1)^j\over 2^j} \binom{k-r}{j}
(r-j)^{2n} .\eqno(2.10)$$
Finally,   the following combinatorial identities hold}
$$\sum_{r=0}^{k} {(-1)^r\over 2^r} \binom{k}{r}\sum_{j=0}^{k-r} {(-1)^j\over 2^j} \binom{k-r}{j}
(r-j)^{2n} = 0, \quad   k= n+1, \ n+2, \dots, \ 2n,\eqno(2.11)$$
$$\sum_{r=0}^{n} {(-1)^r\over 2^r} \binom{n}{r}\sum_{j=0}^{n-r} {(-1)^j\over 2^j} \binom{n-r}{j}
(r-j)^{2n} = (-1)^n n! \  (2n-1)!! \  . \eqno(2.12)$$

\begin{proof} Taking the integral representation (2.2) of $p_n(x)$
and employing the formula (see [6])
$$K_{i\tau}(x)= {\pi\over 2 \sin(\pi i
\tau)}\left[I_{-i\tau}(x)-I_{i\tau}(x)\right],$$
where $I_\nu(z)$ is the modified Bessel function of the third kind
with its series representation
$$I_{-i\tau}(x)= \sum_{k=0}^\infty \frac{(x/2)^{2k-i\tau}}{k!
\Gamma(k+1-i\tau)},\eqno(2.13)$$
where $\Gamma(z)$ is Euler's gamma-function [1], we get the chain of
equalities
$$p_n(x)=  {e^{x}\over 2}\int_{-\infty}^\infty
\frac{(i\tau)^{2n}}{\sin(\pi i \tau)}
\left[I_{-i\tau}(x)-I_{i\tau}(x)\right]\ d\tau =
e^{x}\int_{-\infty}^\infty \frac{(i\tau)^{2n}}{\sin(\pi i
\tau)}I_{-i\tau}(x)\ d\tau
$$
$$= {e^{x}\over \pi i }\int_{-i\infty}^{i\infty} s^{2n-1} \Gamma(1+s)\Gamma(1-s)
\sum_{k=0}^\infty \frac{(x/2)^{2k- s}}{k! \Gamma(k+1-s)}\ ds.
$$
The change of the order of integration and summation in the
right-hand side of the latter equality is indeed possible by
Fubini's theorem owing to the absolute convergence for each $x
>0$, namely
$$\int_{-i\infty}^{i\infty} \left|s^{2n-1}
\Gamma(1+s)\Gamma(1-s)\right| \sum_{k=0}^\infty
\left|\frac{(x/2)^{2k- s}}{k! \Gamma(k+1-s)}\ ds\right| < \infty.$$
Therefore,
$$p_n(x)= {e^{x}\over \pi i}\sum_{k=0}^\infty \frac{(x/2)^{2k}}{k!}\int_{-i\infty}^{i\infty} s^{2n-1}
\frac{\Gamma(1+s)\Gamma(1-s)}{\Gamma(k+1-s)}\left({x\over
2}\right)^{-s}ds.\eqno(2.14)
$$
But the integral in (2.14) can be calculated via the Slater theorem
(see, for instance, in [8], Vol. III) after $2n-1$ times
differentiation with respect to $z$ under the integral sign. This
operation is allowed by virtue of the absolute and uniform
convergence.  Thus we obtain
$$\int_{-i\infty}^{i\infty} s^{2n-1}
\frac{\Gamma(1+s)\Gamma(1-s)}{\Gamma(k+1-s)}z^{-s}ds = - \left(z
{d\over dz}\right)^{(2n-1)}\int_{-i\infty}^{i\infty}
\frac{\Gamma(1+s)\Gamma(1-s)}{\Gamma(k+1-s)}z^{-s}ds$$
$$= 2\pi i \left(z
{d\over dz}\right)^{(2n-1)}\sum_{m=0}^\infty
\frac{(-1)^{m+1}z^{m+1}\Gamma(2+m)}{m!\ \Gamma(2+k+m)}$$
$$= 2\pi i \left(z
{d\over dz}\right)^{(2n-1)}\sum_{m=1}^\infty \frac{(-1)^{m}z^{m}\
m}{(m+k)!}.$$
Hence differentiating under the series sign due to the absolute and
uniform convergence, letting $z= x/2$ and combining with (2.13),
(2.14), we derive the representation
$$p_n(x)= 2e^{x}\sum_{k=0}^\infty \frac{(x/2)^{2k}}{k!}\sum_{m=1}^\infty \frac{(-1)^{m}(x/2)^{m}\
m^{2n}}{(m+k)!} = 2e^{x}\sum_{m=1}^\infty
(-1)^{m}m^{2n}I_m(x),\eqno(2.15)$$
where the inversion of the summation order is motivated by the
absolute convergence of the iterated series. Thus we proved (2.8).

On the other hand, calling relation (5.8.5.3) in [8], Vol. II, it is not difficult to find via differentiation with respect
 to a parameter owing to the absolute and uniform convergence that
$$2 \sum_{m=1}^\infty  m^{2n}I_m(x) = \lim_{a \to 0}  {d^{2n}\over da^{2n}} \ e^{x\cosh a}.$$
We calculate the $2n$-th derivative in the right-hand side of the latter equality employing the Hoppe formula [5]. Precisely it gives,
$$ {d^{2n}\over da^{2n}}  e^{x\cosh a}  =  e^{x\cosh a} \sum_{k=0}^{2n}  x^k \sum_{j=0}^{k}
 {(-\cosh a)^{k-j}\over j! (k-j)!} {d^{2n}\over da^{2n}} \cosh^j a$$
 $$=e^{x\cosh a} \sum_{k=0}^{2n}  x^k \sum_{j=0}^{k}
 {(-\cosh a)^{k-j}\over  2^j j! (k-j)!} {d^{2n}\over da^{2n}}  \sum_{r=0}^{j}  \binom{j}{r} e^{(2r-j)a} $$
$$=e^{x\cosh a} \sum_{k=0}^{2n}  x^k \sum_{j=0}^{k}
 {(-\cosh a)^{k-j}\over  2^j j! (k-j)!}   \sum_{r=0}^{j}  \binom{j}{r} (2r)^{2n-m} (2r-j)^{2n} e^{(2r-j)a}.$$
Therefore, passing to the limit when $a \to 0$, we find the identity
$$2 \sum_{m=1}^\infty  m^{2n}I_m(x) =
 =   e^x \sum_{k=0}^{2n}  {(-x)^k\over k!} \sum_{j=0}^{k}  \binom{k}{j} {(-1)^{j} \over  2^j }
 \sum_{r=0}^{j}   \binom{j}{r}  (2r -j)^{2n}$$
$$=   e^x \sum_{k=0}^{2n}  {(-x)^k\over k!}\sum_{r=0}^{k} {(-1)^r\over 2^r} \binom{k}{r}\sum_{j=0}^{k-r} {(-1)^j\over 2^j} \binom{k-r}{j} (r-j)^{2n}.$$
On the other hand, since $(-1)^{m}I_m(x)= I_m(-x)$, we return to (2.15) and appealing  to the right-hand side of the latter equality, we derive
$$p_n(x)=   \sum_{k=0}^{2n}  {(-x)^k\over k!}\sum_{r=0}^{k} {(-1)^r\over 2^r} \binom{k}{r}\sum_{j=0}^{k-r} {(-1)^j\over 2^j} \binom{k-r}{j} (r-j)^{2n}.$$
However, as it is proved in [11],  $p_n(x)$ is a polynomial of degree $n$.  Therefore,  all coefficients in front of powers $x^{k}, \ k= n+1,\  n+2, \  \dots, \  2n$ are  surely equal to zero.   Thus we establish the explicit formula (2.9) with coefficients (2.10) and combinatorial equality (2.11).  Taking into account the value (2.4) of the leading coefficient $a_{n,n}$, we get identity (2.12) and complete the proof of Theorem 1.
\end{proof}

Finally, in this section we will established an exact asymptotic behavior of $p_n(x)$, when $n \to \infty$ and $x >0$ is a fixed number.    Precisely, we have

{\bf Theorem 2.}  {\it Let $x >0$ and $\beta \in (0, \pi/2)$ be fixed numbers. Then}
$$p_n(x)=    {6x(-1)^n \sin\beta \  (2n)!\over \pi \beta^{2n} (2n+1)^3} e^{x} \left( 1+ O\left({1\over n}\right)\right), \ n \to \infty.\eqno(2.16)$$

\begin{proof}  Indeed, employing representation (2.2) and choosing a fixed parameter $\beta \in (0, \pi/2)$, we write
$$p_n(x)= {2(-1)^n\over \pi} e^{x}\int_{0}^\infty \tau^{2n} e^{-\beta\tau} \left[\cosh(\beta\tau) + \sinh(\beta\tau)\right]
K_{i\tau}(x)\ d\tau. $$
Hence the Parseval identities  for the cosine and sine  Fourier
transforms [9] together with relations (2.5.31.4) in [8], Vol. I and
(2.16.48.20) in [8], Vol. II drive us to the equalities
$$p_n(x)= {2(-1)^n (2n)!\over \pi} e^{x}\int_{0}^\infty {e^{-x\cos(\beta )\cosh
y}\over (\beta^2+
y^2)^{n+1/2}}\left[\cos\left((2n+1)\hbox{arctg}(y/\beta)\right)\cos(
x\sin(\beta)\sinh y )\right.$$
$$\left. + \sin\left((2n+1)\hbox{arctg}(y/\beta)\right)\sin(
x\sin(\beta)\sinh y )\right]dy$$
$$= {2(-1)^n (2n)!\over \pi} e^{x}\int_{0}^\infty {e^{-x\cos(\beta )\cosh
y}\over (\beta^2+ y^2)^{n+1/2}}\cos\left[(2n+1)\hbox{arctg}(y/\beta)
- x\sin(\beta)\sinh y\right]dy$$
$$= {2 (-1)^n (2n)!\over \pi} e^{x}\  {\rm Re} \left[ \int_{0}^\infty e^{-x\cos(\beta+iy)} \exp \left[-(2n+1)\left(\log\left(\sqrt{\beta^2+ y^2} \right)\right.\right.\right.$$
$$\left.\left.\left.  + i\hbox{arctg}(y/\beta)\right)\right]dy\right] .\eqno(2.17)$$
The asymptotic behavior for large $n$ of the latter integral in brackets can be treated by the Laplace method [7], Ch. 4.
Indeed, it has the form
$$\int_{0}^\infty q(y) e^{-(2n+1)p(y)} dy,$$
where $p(y)=  \log(\beta +iy), \ q(y)= e^{-x\cos(\beta+iy)}$.   Thus the integral
$$\int_{0}^\infty {e^{-x\cos(\beta+iy)} \over (\beta+ iy)^{2n+1}} dy \sim {1\over \beta^{2n+1}}
\sum_{s=0}^\infty {a_s(x,\beta)   \ s!\over (2n+1)^{s+1}},\  n \to \infty,$$
where the calculation of three first coefficients $a_s$ will be enough to determine the main term of the expansion of its real part.  Namely, using formulas for coefficients in [7], we have
$$a_0=  -i\beta e^{-x\cos \beta}, \quad a_1= - i\beta (1 + x\beta \sin\beta)  \  e^{-x\cos \beta},$$
$$a_2= {3\over 2}x\beta\sin\beta - {i \beta\over 2} e^{-x\cos
\beta}\left(\beta^2 x(\cos\beta + x\sin^2\beta) + 1\right).$$
Therefore,   returning  to (2.17),  we  find

$$p_n(x)=  {6x (-1)^n \sin\beta \  (2n)!\over \beta^{2n} (2n+1)^3 \pi} e^{x}\left( 1+ O\left({1\over n}\right)\right), \ n \to \infty,$$
which proves (2.17).

\end{proof}

\section {Representations of the Yor integral}

This section will complete our goal to represent integral (1.1) in a
different form. Namely, we will start relating $F_t(r)$ with the
inverse Kontorovich-Lebedev transform (1.4). Indeed, taking (1.4)
with $\nu=i\tau$ and integrating by parts we come out with the
representation of the modified Bessel function
$$\tau K_{i\tau}(r)=r\int_0^\infty e^{-r\cosh u}\sinh u\sin\tau u du, \ r>0. \eqno(3.1)$$
Hence applying the Parseval equality for the sine Fourier transform
[9] with the relation (2.5.36.1) in [8], Vol. 1, integral (1.1)
becomes $(r,t >0)$
$$F_t(r)= \frac{  e^{\pi^2/2t}}{r \pi^2}\int_0^\infty \tau K_{i\tau}(r) \left[ \exp\left(- {t
(\tau- \pi/t)^2\over 2}\right) - \exp\left(- {t (\tau +
\pi/t)^2\over 2}\right)\right]d\tau$$
$$=\frac{ 2}{r \pi^2}\int_0^\infty e^{- {t\over 2} \tau^2}\tau \sinh\pi\tau K_{i\tau}(r)d\tau.\eqno(3.2)$$
Hence form (1.2) we obtain reciprocally for each $t >0$
$$\int_0^\infty K_{i\tau}(r)F_t(r) dr = e^{- {t\over 2} \tau^2}\eqno(3.3)$$
and via Parseval equality (1.3) we derive the value of the integral
$$\int_0^\infty |F_t(r)|^2 r dr  = {2\over \pi^2}\int_0^\infty \tau\sinh(\pi\tau)\  e^{- t \tau^2} d\tau$$
$$= {2 e^{\pi^2/4t} \over \pi^2}\int_{-\infty}^\infty \tau  e^{- t (\tau- \pi/(2t))^2}
d\tau = {e^{\pi^2/4t} \over  t\sqrt {\pi t}}.$$

Recently (see [12]), the author introduced the following heat kernel
for the Kontorovich-Lebedev transform
$$ h_t(x,y)= {2\over x\pi^2}\int_0^\infty e^{-t\tau^2/2}\tau
\sinh\pi\tau \ K_{i\tau}(x)K_{i\tau}(y)d\tau.\eqno(3.4)$$
Hence appealing to the Macdonald formula (1.15) and Fubini's theorem
to interchange the order of integration,  we come out with the
representation of (3.4) as a translation operator of the Yor
integral for convolution (1.12). Precisely, minding (3.2) it gives
$$h_t(x,y)= {1\over 2x}\int_0^\infty e^{-{1\over
2}\left(r{x^2+y^2\over xy}+{xy\over r}\right)}F_t(r) dr ,\
x,y >0.\eqno(3.5)$$

Meanwhile differential and convolution properties of the Yor integral (1.1) are
given by the following

{\bf Theorem 3}. {\it The function $F_t(r)$ is infinitely
differentiable  of  variables $(r,t) \in \mathbb{R}_+\times
\mathbb{R}_+$, satisfying the estimate
$$\left|{\partial^m F_t(r)\over \partial t^m}\right| \le
{2^{1/4-m}e^{\pi^2/t}\over t^{m+ 3/4}\pi^{11/8}}\
K_0^{1/2}(2r)\Gamma^{1/4}(4m+5/2), \quad m \in
\mathbb{N}_0.\eqno(3.6)$$
Moreover,  $F_t(r)$ is a solution of the generalized diffusion
equation $(u= u(t,r))$
$$2{\partial u\over \partial t}= r^2{\partial^2 u\over \partial r^2} + r{\partial u\over \partial r}
- r^2 u \eqno(3.7)$$
and satisfies the index law in terms of convolution $(1.12)$
$$F_t(r)= \left(F_{t_1}* F_{t_2}\right)(r), \quad t_1+t_2=t .\eqno(3.8)$$
Finally,  for the Yor integral  the following integral equation takes place}
$$\int_0^\infty h_{t_1}(r,y)F_{t_2}(y)dy= F_{t_1+ t_2}(r). \eqno(3.9)$$

\begin{proof}  In fact, employing the following inequality (see [11], [13]) for
derivatives of the modified Bessel function with respect to $x$
$$\left|{\partial^m K_{i\tau}(r)\over \partial x^m} \right| \le \ e^{-\delta \tau}K_m(r \cos\delta), \ x >0, \ \tau
>0, \ \delta \in \left[0; {\pi\over 2}\right), \ m= 0,1,\dots \eqno(3.10)$$
it is not difficult to verify that all positive $t$ the integral in
the right-hand side of the latter equality in (3.2) together with
derivatives of any order with respect to $r$ converge absolutely and
uniformly by $r \ge r_0 > 0$. Therefore $F_t(r)$  is infinitely
differentiable with respect to $r$. Similar motivation can be done
for the derivatives of the order $m \in \mathbb{N}_0$ with respect
to $t >0$ and it gives the expression
$${\partial^m F_t(r)\over \partial t^m} = {2^{1-m} (-1)^m \over \pi^2}\int_0^\infty
e^{-t\tau^2/2}\tau^{2m+1} \sinh\pi\tau \ K_{i\tau}(r)d\tau, \ m=
0,1,\dots\ .\eqno(3.11)$$
Hence the Schwarz inequality and relation (2.16.52.6) in [8], Vol.
II yield
$$\left|{\partial^m F_t(r)\over \partial t^m}\right| \le {2^{1-m}\over \pi^2}\left(\int_0^\infty
e^{- t\tau^2+ 2\pi\tau}\tau^{2(2m+1)} d\tau\right)^{1/2}
\left(\int_0^\infty K_{i\tau}^2(r)d\tau\right)^{1/2}$$
$$\le {2^{1/2-m}\over \pi\sqrt\pi}K_0^{1/2}(2r) \left(\int_{-\infty}^\infty
e^{- t\tau^2+ 4\pi\tau} d\tau\right)^{1/4}\left(\int_0^\infty e^{-
t\tau^2}\tau^{4(2m+1)} d\tau\right)^{1/4}$$
$$={2^{1/4-m}\over t^{m+ 3/4}\pi^{11/8}}\ e^{\pi^2/t}
K_0^{1/2}(2r)\Gamma^{1/4}(4m+5/2).$$
Thus we proved  (3.6).

Further, as it follows from (1.5), (1.6) and absolute and uniform
convergence of the corresponding integrals, formula (3.11) can be
written as the following partial differential equation
$${\partial^m F_t(r)\over \partial t^m} = (-1)^m 2^{-m} {\cal A}_r^m F_t(r),\eqno(3.12)$$
where $\ m= 0,1,\dots\ $ and $A_r^m$ is $m$-th iterates of the
operator (1.5). In particular, letting $m=1$ we obtain that the Yor
integral (1.1) satisfies the generalized diffusion equation (3.7).

Next, equality (3.8) is a direct consequence of (3.3), factorization equality (1.14) and the uniqueness property of the Kontorovich-Lebedev transform [11].  In order to prove (3.9),  we  multiply  $h_{t_1}(x,y)$ by $F_{t_2}(y)$ and integrate by $y$ over $\mathbb{R}_+$. Then  using   (3.3), (3.4) and the Fubini theorem, which is applicable due to the absolute convergence of the iterated integral (it can be verified,  employing inequalities (3.6) with $m=0$ and (3.10)), we get the result.
\end{proof}

Finally, we will establish a representation of  the Yor integral (1.1) in terms of polynomials (2.1) and the heat kernel (3.4).

We have

{\bf Theorem 4}. {\it Let $r, t > 0$. Then the Yor integral satisfies the following  equation
$$F_t(r) =  \sum_{k=1}^\infty {(-1)^k  \pi^{2(k-1)} \over (2k-1)!}  \left(e^{-r} {p_k(r)\over r} * F_t(r)\right),\eqno(3.14)$$
where $*$ denotes convolution $(1.12)$.  Moreover, it can be rewritten in the form
$$ F_t(r) =   \sum_{k=1}^\infty {(-1)^k  \pi^{2(k-1)}a_k(r,t) \over (2k-1)!},\eqno(3.15)$$
where
$$a_k(r,t)= \int_0^\infty e^{-u} h_t(r,u)  p_k(u){du\over u}$$
and $h_t$ is the heat kernel $(3.4)$.}

\begin{proof}  Calling again integral (3.2),   we  expand the hyperbolic sine  in Taylor's series. Then changing the order of integration and summation due to the estimate
$$\int_{0}^\infty e^{-t\tau^2/2} \tau |K_{i\tau}(r)|  \sum_{k=1}^\infty {(\pi\tau)^{2k-1} \over (2k-1)!} d\tau
\le K_0(r) \sum_{k=1}^\infty {\pi^{2k-1} \over (2k-1)!} \int_{0}^\infty e^{-t\tau^2/2} \tau^{2k} d\tau$$
$$= {K_0(r)\over \pi \sqrt {2t} } \sum_{k=1}^\infty {\left(\pi\sqrt{2/t}\right) ^{2k} \Gamma(k+ 1/2) \over (2k-1)!} < \infty, $$
we come out with the representation
$$F_t(r)= \frac{2}{r \pi^2}   \sum_{k=1}^\infty {\pi^{2k-1} \over (2k-1)!} \int_{0}^\infty  e^{-t\tau^2/2}
\  \tau^{2k}  K_{i\tau}(r)d\tau.\eqno(3.16)$$
But the integral in (3.16) can be treated  with the convolution for the Kontorovich-Lebedev transform (1.12),
its factorization property (1.14)  and  inversion formula (1.4). Thus minding formulas (2.2),  (3.2) and (3.3), we derive
$$ \frac{2}{r \pi^2} \int_{0}^\infty  e^{-t\tau^2/2} \  \tau^{2k}  K_{i\tau}(r)d\tau=  {(-1)^k\over \pi}
\left(e^{-r} {p_k(r)\over r} * F_t(r)\right). $$ Substituting in
(3.16), it gives
$$F_t(r) =   \sum_{k=1}^\infty {(-1)^k  \pi^{2(k-1)} \over (2k-1)!}  \left(e^{-r} {p_k(r)\over r} * F_t(r)\right),$$
and equality (3.14) is proved.  On the other hand,  employing (3.5),
one can express the latter convolution in terms of the heat kernel (3.4).  Precisely, we find  for each $k \in \mathbb{N} $ (see (1.12))
$$\left(e^{-r} {p_k(r)\over r} * F_t(r)\right) = \int_0^\infty e^{-u} h_t(r,u)  p_k(u){du\over u}= a_k(t,r),$$
which proves  (3.15).

\end{proof}

\bigskip
\centerline{{\bf Acknowledgments}}
\bigskip
The present investigation was supported, in part,  by the "Centro de
Matem{\'a}tica" of the University of Porto.

\bigskip
\centerline{{\bf References}}
\bigskip
\baselineskip=12pt
\medskip
\begin{enumerate}

\item[{\bf 1.}\ ]
 M. Abramowitz, I.A. Stegun,
{\it Handbook of Mathematical Functions}, Dover, New York, 1972.

\item[{\bf 2.}\ ] H. Alzer,  Sharp bounds for the Bernoulli numbers, {\it Arch. Math.} {\bf
74} (2000),  3,  207-211.

\item[{\bf 3.}\ ] P. Barrieu, A. Rouault and M. Yor, A study of Hartman-Watson distribution motivated by numerical problems
related to Asian options pricing, {\it Journ. Appl. Probab.} {\bf
41} (2004), 7,  1049-1058.

\item[{\bf 4.}\ ] P. Carr and M. Schr$\ddot{o}$der, Bessel processes, the integral and geometrical Brownian motion, and Asian options, {\it Theory Probab. Appl.} {\bf 48} (2004), 3, 400-425.

\item[{\bf 5.}\ ]   W.P. Johnson,   The curious history of Faa di Bruno's formula,  {\it The Amer. Math. Monthly,} {\bf 109}  (2002), 3,  217- 234.

\item[{\bf 6.}\ ]
 N.N. Lebedev, {\it Special Functions and Their Applications}, Dover,  New York, 1972.

 \item[{\bf 7.}\ ]
 F.W.J. Olver, {\it Introduction to Asymptotics and Special Functions},  Academic Press,  New York and London, 1974.

\item[{\bf 8.}\ ] A.P. Prudnikov, Yu.A. Brychkov and O.I.
Marichev, {\it Integrals and Series}. Vol. I: {\it Elementary
Functions}, Vol. II: {\it Special Functions}, Gordon and Breach, New
York and London, 1986, Vol. III: {\it More Special Functions},
Gordon and Breach, New York and London, 1990.

\item[{\bf 9.}\ ]
 I.N. Sneddon, {\it The Use of Integral  Transforms}, McGraw-Hill,
 New York, 1972.

\item[{\bf 10.}\ ]S.B. Yakubovich and Yu.F. Luchko, {\it The
Hypergeometric Approach to Integral Transforms and Convolutions},
(Kluwers Ser. Math. and Appl.: Vol. 287), Dordrecht, Boston, London,
1994.

\item[{\bf 11.}\ ] S.B. Yakubovich, {\it Index Transforms}, World
Scientific Publishing Company, Singapore, New Jersey, London and
Hong Kong, 1996.

\item[{\bf 12.}\ ] S.B. Yakubovich, A class of polynomials and discrete transformations associated with the
Kontorovich-Lebedev operators, {\it Integral Transforms and Special
Functions} {\bf 20} (2009), 7,  551- 567.

\item[{\bf 13.}\ ] S. Yakubovich, The heat kernel and Heisenberg inequalities related to the Kontorovich-Lebedev transform,  {\it Commun. Pure Appl. Anal.} {\bf 10} (2011), N 2, 745-760.

\item[{\bf 14.}\ ] M. Yor, Loi de l'indice du lacet Brownien et distribution de Hartman-Watson, {\it Z. Wahrscheinlichkeitsth. } {\bf 53} (1980), 71-95 (in French).

\item[{\bf 15.}\ ] M. Yor, On some exponential functionals of Brownian motion, {\it Adv. Appl. Probab., } {\bf
24} (1992), 509-531.

\end{enumerate}

\vspace{5mm}

\noindent Semyon  Yakubovich\\
Department of  Mathematics,\\
Faculty of Sciences,\\
University of Porto,\\
Campo Alegre st., 687\\
4169-007 Porto\\
Portugal\\
E-Mail: syakubov@fc.up.pt\\

\end{document}